\documentclass[10pt]{amsart}
\usepackage{amsmath,amssymb,amsfonts} % Typical maths resource packages
\usepackage{graphics}                 % Packages to allow inclusion of graphics
\usepackage[square, numbers]{natbib}
\usepackage[english]{babel}
\usepackage{tikz} 
\usepackage[utf8]{inputenc}
\usepackage[all]{xy}	
\usepackage[T1]{fontenc}
\usepackage{pifont}

\usepackage[margin=1.2in]{geometry}

\newcommand{\C}{\mathbb{C}}

\newcommand{\R}{\mathbb{R}}

\theoremstyle{definition}
\newtheorem{defn}{Definition}[subsection]
\newtheorem{rem}[defn]{Remark}
\newtheorem{cons}[defn]{Construction}
\newtheorem{exam}[defn]{Example}

\theoremstyle{plain}
\newtheorem{Prop}[defn]{Proposition}
\newtheorem{Thm}[defn]{Theorem}
\newtheorem{cor}[defn]{Corollary}
\newtheorem{lem}[defn]{Lemma}

\newtheorem*{asser}{Claim}

\theoremstyle{exampstyle}\newtheorem{thmx}{Theorem}

\newtheorem{Propx}{Proposition}

\theoremstyle{exampstyle}

\newtheorem{corx}{Corollary}

\def\Ind{\setbox0=\hbox{$x$}\kern\wd0\hbox to 0pt{\hss$\mid$\hss}
\lower.9\ht0\hbox to 0pt{\hss$\smile$\hss}\kern\wd0}
\def\Notind{\setbox0=\hbox{$x$}\kern\wd0\hbox to 0pt{\mathchardef
\nn=12854\hss$\nn$\kern1.4\wd0\hss}\hbox to
0pt{\hss$\mid$\hss}\lower.9\ht0 \hbox to
0pt{\hss$\smile$\hss}\kern\wd0}

\newtheoremstyle{exampstyle}
{3pt} % Space above
{3pt} % Space below
{\itshape} % Body font
{} % Indent amount
{\bfseries} % Theorem head font
{.} % Punctuation after theorem head
{.5em} % Space after theorem head
{} % Theorem head spec (can be left empty, meaning `normal')

\makeindex

\begin{document}

\title[Rational factors, invariant foliations and Anosov flows]{Rational factors, invariant foliations and algebraic disintegration of compact mixing Anosov flow of dimension $3$}
\author{Rémi Jaoui}
\address{Rémi Jaoui, Pure Mathematics Department, University of Waterloo, Ontario, Canada}
\email{rjaoui@uwaterloo.ca}
\date\today

\maketitle

\begin{abstract}
In this article, we develop a geometric framework to study the notion of semi-minimality for the generic type of a smooth autonomous differential equation $(X,v)$, based on the study of rational factors of $(X,v)$ and of algebraic foliations on $X$, invariant under the Lie-derivative of the vector field $v$.

We then illustrate the effectiveness of these methods by showing that certain  autonomous algebraic differential equation of order three defined over the field of real numbers --- more precisely, those associated to mixing, compact, Anosov flows of dimension three --- are generically disintegrated.   
\end{abstract}

\vspace{0.2cm}

Many applications of geometric stability theory to concrete geometric settings originated from Zilber's trichotomy on minimal types --- respectively called disintegrated, locally-modular (non-disintegrated) or non-locally modular --- in stable theories. From this point of view, the study of ordinary algebraic differential equations takes place in the theory  $\textbf{DCF}_0$ of existentially closed differential fields of characteristic $0$, where the minimal types of the second and of the third kinds have been entirely classified in \cite{Sok}.

This theorem entails many important consequences for minimal algebraic differential equations while certain algebraic differential equations of higher rank witness a ``mixed'' behavior where pregeometries of different kinds are witnessed. In the still unpublished \cite{moi}, we advocated that many interesting consequences of this theorem remain valid under a weaker minimality notion, sufficient to rule out these mixed behaviors, known as \textit{semi-minimality} (see also Section 2.5 of \cite{moi4} for a geometric exposition).

In this article, we develop a complementary geometric framework for the notion of the semi-minimality, which provide effective tools to establish this property, based on the study of \textit{rational factors and invariant foliations  of a given autonomous differential  equations $(X,v)$}.

We also apply these techniques to study certain three dimensional algebraic differential equations --- namely, \textit{algebraically presented, compact, mixing Anosov flows of dimension three} ---  and establish \textit{disintegration} of the system of algebraic relations shared by their generic solutions. The techniques of this article will also be applied  in \cite{moi4} to study similar properties for generic planar algebraic vector fields \\

\textbf{Semi-minimality and rational factors.}  We fix $k$ a field of characteristic $0$ and we consider autonomous algebraic differential equations $(X,v)$ with coefficients in $k$, presented as smooth irreducible $k$-algebraic varieties $X$ endowed with algebraic vector fields $v$.

The solutions of such a differential equation $(X,v)$ in an existentially closed differential field form a definable set of the theory $\textbf{DCF}_0$. All the generic solutions of $(X,v)$ over $k$ --- those who do not lie in any proper closed $k$-subvariety of $X$ --- realize the same (complete) type $p \in S(k)$ called \textit{the generic type of $(X,v)$}. Recall also that a type $p \in S(k)$ in a stable theory is called \textit{semi-minimal} if it is almost internal to the collection of $\mathrm{acl}(k)$-conjugates of a minimal type $r \in S(l)$ defined over some extension of the parameters $k \subset l$.  

As many other notions from geometric stability, this definition involves, as it stands, arbitrary extension of parameters, or in other words, base change by arbitrary differential field extensions of $(k,0)$. It is remarkable that, in fact, this property can be checked directly from $k$-algebraic dynamical data (not involving any base-change), using the notion of rational factor as follows.

A \textit{rational factor} $\phi: (X,v) \dashrightarrow (Y,w)$ (over $k$) of an autonomous algebraic differential equation $(X,v)$ is a rational dominant morphism $\phi: X \dashrightarrow Y$ over $k$, satisfying the obvious compatibility relation $d\phi(v) = w$ with the vector fields $v$ and $w$. Additionally, we say that an autonomous differential equation $(X,v)$ over $k$ \textit{does not admit any non-trivial rational factor} if any rational factor of $(X,v)$  of positive dimension has a finite generic fibre.

With this terminology in place, it follows from standard techniques of geometric stability theory that the generic type of an algebraic differential which does not admit any non-trivial rational factor is always semi-minimal (see, for example, \cite{Moosa1}). The first aim of this article is to develop a correspondent geometric framework to study effectively rational factors of algebraic autonomous differential equations. \\
  
%One of the main obstacle to a direct dynamical analysis of the rational factors $\pi: (X,v) \dashrightarrow (Y,w)$ of an autonomous differential equation $(X,v)$ lies in the indeterminacy locus of $\pi$. In this article, we develop and study a notion of \textit{invariant foliation} --- analogous to the classical notion of invariant subvariety --- for smooth differential equation $(X,v)$ in order to overcome this obstacle.\\

\textbf{Invariant foliations.} Our analysis of rational factors of an algebraic autonomous differential equation $(X,v)$ relies on an auxiliary construction which associates to a rational factor of relative dimension $r$, an invariant foliation of $(X,v)$ of rank $r$, that we describe below. \\

Let $X$ be a smooth complex algebraic variety over some field $k$ of characteristic $0$. Recall that a non-singular foliation on $X$ is an involutive vector subbundle of the tangent bundle $T_{X/k}$ on $X$. More generally, a \textit{(possibly singular) foliation $\mathcal F$ of rank $r$} on $X$ is a coherent, integrable, saturated subsheaf of rank $r$ of the locally free sheaf $\Theta_{X/k}$ of vector fields on $X$. 

It follows from the results of \cite{Har2} on coherent saturated subsheaves of a locally free sheaf that a (possibly singular) foliation $\mathcal F$ is always non-singular outside of a closed subvariety $Z = \mathrm{Sing}(\mathcal F)$ of $X$ of codimension $\geq 2$, called the \textit{singular locus of $\mathcal F$}. Our main incentive to allow such singularities in this definition of foliation is the following extension principle for foliations (Proposition \ref{saturation}), which allows us to localize our analysis and constructions at the generic point of $X$: \textit{provided $X$ is smooth, any foliation on a dense open set $U \subset X$ extends uniquely to a foliation  on $X$}. \\

Assume now moreover that the smooth algebraic variety $X$ in endowed with a vector field $v$. Recall that the Lie derivative $\mathcal L_v$ of the vector field $v$ acts on the coherent sheaf $\Theta_{X/k}$ of vector fields  on $X$ by the rule $\mathcal L_v(w) = [v_{|U},w]$ for every local section $w \in \Theta_{X/k}(U)$.

For every open set $U$ of $X$, the Lie-derivative $\mathcal L_v$ acts as a derivation on the $\mathcal O_X(U)$-module $\Theta_{X/k}(U)$, compatible with the derivation $\delta_v$ induced by the vector field $v$ on $\mathcal O_X(U)$.  In other words, the Lie derivative $\mathcal L_v$ defines a  $\delta_v$-module structure (in the sense of Pillay-Ziegler \cite{Pil}) on $\Theta_{X/k}(U)$ viewed as an $(\mathcal O_X(U),\delta_v)$-module. We study, more generally, in the first section of this article coherent sheaves $\mathcal E$ on $X$ endowed with such a connection operator $\nabla$ over an autonomous differential equation $(X,v)$, under the name of $D$-coherent sheaves $(\mathcal E,\nabla)$ over the $D$-scheme $(X,v)$. \\

A foliation $\mathcal F$ on $X$ is called an \textit{invariant foliation} of $(X,v)$ if it is stable under the action of the Lie-derivative $\mathcal L_v$ or, in other words, if $\mathcal F$ is a $D$-coherent subsheaf of the $D$-coherent sheaf $(\Theta_{X/k},\mathcal L_v)$. Among other things, we provide, for smooth complex algebraic autonomous differential equations $(X,v)$, an analytic interpretation (Proposition \ref{analytification invariance}) for this notion of invariance  and we prove that (Corollary \ref{invariance-singular}) \textit{the singular locus of an invariant foliation $\mathcal F$ of $(X,v)$ is always a closed invariant subvariety of $(X,v)$}.

The following proposition relates the study of rational factors of smooth algebraic autonomous differential equations $(X,v)$ with the study of invariant foliations of $(X,v)$.

{\Propx  Let $k$ be a field of characteristic $0$ and let $(X,v)$ be a $k$-algebraic irreducible variety endowed with a vector field $v$.

If $\phi :(X,v)\dashrightarrow (Y,w)$ is a rational factor of $(X,v)$ of relative dimension $r$ then the foliation $\mathcal F_{\phi}$ tangent to the fibres of $\phi$ is an invariant foliation of $(X,v)$ of rank $r$.} 
\vspace{0.2cm}

Recall that, over a field $k$ of characteristic $0$, any dominant rational morphism $\phi: X \dashrightarrow Y$ of relative dimension $r$ determines a foliation $\mathcal F_\phi$ on $X$ of rank $r$, called the \textit{foliation tangent to the fibre of $\phi$}: indeed, one can first define $\mathcal F_{\phi |U}$ on any dense open set $U \subset X$ on which $\pi$ is defined and smooth and then extend it to $X$ itself using the extension principle discussed above. Note also that, in fact, this construction presents the foliation $\mathcal F _\phi$ of rank $r$ as a fibration of $X$ by algebraic subvarieties of dimension $k$. \\

\textbf{Disintegration of $3$-dimensional compact mixing Anosov flow.} The second aim of this article is to exploit the relationship between invariant foliations and rational factors described above, conjointly with the results of \cite{moi2}, to study disintegration properties of, algebraically presented, compact, mixing Anosov flows of dimension $3$. \\

Recall first the notion of disintegration (as formulated, for example, in \cite{moi}) at the generic point of an autonomous differential equation $(X,v)$ (over some field $k$ of characteristic $0$). For $n \geq 2$, denote by $\mathcal I^{gen}_n(X,v)$ the set of closed irreducible invariant subvarieties of $(X,v)^n$  which project generically on all factors. We say that an  algebraic differential equation $(X,v)$ is \textit{generically disintegrated} if for every $n \geq 3$, every $Z \in \mathcal I^{gen}_n(X,v)$ can be written as an irreducible component (which projects generically on all factors) of: 
$$\bigcap_{1 \leq i \neq j \leq n} \pi_{i,j}^{-1}(Z_{i,j}).$$
where $\pi_{i,j} : X^n \longrightarrow X^2$ is the projection on the $i^{th}$ and $j^{th}$ coordinates and $Z_{i,j} \in \mathcal I^{gen}_2(X,v)$ for every $i \neq j$.

The following theorem shows that an algebraically presented, mixing, compact Anosov flow of dimension three is always generically disintegrated:

{\thmx\label{maintheoremInt} Let $X$ be an absolutely irreducible variety of dimension $3$ over $\mathbb{R}$ endowed a vector field $v$. Assume that the real-analytification $X(\R)^{an}$ of $X$ admits a compact (non-empty) connected component $C_\mathbb{R}$ contained in the regular locus of $X$.

If the real-analytic flow $(C_\mathbb{R}, (\phi_t)_{t \in \R})$ is a mixing Anosov flow, then the autonomous differential equation $(X,v)$ is generically disintegrated.}
\vspace{0.3cm}

Anosov flows are uniformly hyperbolic flows defined by Anosov in \cite{Ano}. Recall that, for an Anosov flow $(M,(\phi_t)_{t \in \R})$, the various notions of mixing --- topologically weakly-mixing, topologically mixing and the mixing properties relatively to an equilibrium measure--- collapse into a single one (see for example \cite{Coud}). We will simply say that  $(M,(\phi_t)_{t \in \R})$ is a \textit{mixing Anosov flow} to mean that one of the previous properties is satisfied.

Apart from suspensions of Anosov diffeomorphisms, other classical examples of compact Anosov flow come from geodesic motion on compact Riemannian manifolds with negative curvature and are known to be mixing (see \cite{Dalbo}). We proved in \cite{moi2} that these classical examples ensure the existence of unlimited families of algebraic autonomous differential equations satisfying the hypotheses of Theorem \ref{maintheoremInt}. \\

 Theorem \ref{maintheoremInt} provides a substential improvement of some of results of \cite{moi2} which only ensures, for an autonomous differential equation $(X,v)$ satisfying the hypotheses of Theorem \ref{maintheoremInt}, the existence of a generically disintegrated rational factor $\pi:(X,v)\dashrightarrow (Y,w)$ of positive dimension.
 
Building on this result, Theorem \ref{maincorollaryInt} reduces to proving that  an an algebraically presented, compact, mixing Anosov flow $(X,v)$ of dimension three   does not admit any non-trivial rational factor. For that matter, we carry out, in the fourth section, an analysis of the continuous invariant distributions by an Anosov flow $(M,(\phi_t)_{t \in \R})$ of dimension three based on the usual (hyperbolic) decomposition:

$$ TM = W^{ss} \oplus E \oplus W^{su}$$
of the tangent bundle of $M$ as the direct sum of the strongly stable continuous line bundle $W^{ss}$, the strongly unstable continuous line bundle $W^{su}$ and the direction $E$ of flow. 

The decisive argument in our study of rational factors through invariant foliations is a result of Plante in \cite{Plante}, which ensures that the strongly stable and the strongly unstable  distributions  $W^{ss}$ and $W^{su}$  do not admit any algebraic leaves.  \\

To conclude this introduction, we illustrate more concretely the content of Theorem \ref{maintheoremInt} with a separate real-analytic instance of it, when applied to a geodesic motion in negative curvature:

{\corx\label{maincorollaryInt} Let $M$ be a regular compact real-algebraic subset of the Euclidean space $\mathbb{R}^N$ of dimension $2$ with negative curvature and let $r$ be integer $ \geq 2$.

Consider $r$ unitary geodesics $\gamma_1, \ldots \gamma_r: \mathbb{R} \rightarrow SM$ of  the Euclidean submanifold $M$, viewed as analytic curves on the sphere bundle $SM \subset T \mathbb{R}^N$ of $M$. Assume that, for every $1 \leq i \leq r$, the analytic curve $\gamma_i$ is Zariski-dense in $SM$. Then, the following are equivalent:
\begin{itemize}
\item[(i)] The analytic curve $t \mapsto (\gamma_1(t), \ldots, \gamma_r(t))$ is Zariski-dense in $SM^r$.
\item[(ii)] For every $i \neq j$, the analytic curve $t \mapsto (\gamma_i(t),\gamma_j(t))$ is Zariski-dense in $SM^2$.  
\end{itemize}}
\vspace{0.3cm}

\textbf{Organization of the article.} In the first section, we recall various classical facts about the Lie-derivative in the setting of algebraic geometry. We emphasize the connection with differential algebra by formulating the results at the level of $D$-schemes.

In the second section, we set up some classical definitions and notations from the theory of algebraic foliations before defining, in the third section, the notion of invariant algebraic foliation for a vector field $v$ on a smooth algebraic variety $X$. 

Finally the last section is dedicated to the study of rational factors of compact, mixing, Anosov flows of dimension $3$ and to the proofs of Theorem \ref{maintheoremInt} and Corollary \ref{maincorollaryInt}. \\

\textbf{Acknowledgments.} The results of this article are based on the many interesting ideas and suggestions of my PhD supervisors, Jean Benoît Bost and Martin Hils. I would also like to thank them for their numerous comments on earlier versions of this article. I am also very grateful to Rahim Moosa for many useful discussions and comments on the content of this article.

\tableofcontents

\section{Lie derivative}
In this section, we recall some classical facts on the Lie derivative in the setting of algebraic and analytic geometry.

The Lie derivative of a vector field (or more generally of a tensor field) with respect to another vector field is a standard notion of differential calculus. Surprisingly, it seems that that this notion has not been considered before in the setting of differential algebra \textit{à la Buium} \cite{Bui}.

On the other hand, the effectiveness of the Lie derivative in order to study \textit{``non-integrability properties''} of a differential equation given by a vector field on a manifold already appears in the work of Morales-Ruiz and Ramis (cf. \cite{Morales} and \cite{MoralesRamis}) by means of the so-called variational equation (see also \citep[Part 3.1, p.46]{Aud}). Moreover, the notion of Lie-derivative lies --- through Frobenius Integrability Theorem --- at the heart of the theory of algebraic foliations of dimension $\geq 2$ , that we will study in the second section of this article.

\subsection{Definition} We fix $k$ a field of characteristic $0$. Recall that if $A$ is a $k$-algebra then the $A$-module $\mathrm{Der}_k(A) = \mathrm{Hom}(\Omega^1_{A/k};A)$ is endowed with a Lie-bracket by the formula:

\begin{eqnarray}\label{Liebracket}
[\delta_1, \delta_2] = \delta_1 \circ \delta_2 - \delta_2 \circ \delta_1.
\end{eqnarray}

The Lie-bracket is compatible with localization. Indeed, if $A$ is a $k$-algebra and $S$ is a multiplicative system of $A$ then the natural isomorphism of $S^{-1}A$-modules:
$$ S^{-1}\mathrm{Der}_k(A) \simeq \mathrm{Der}_k(S^{-1}A)$$
is also an isomorphism of $k$-Lie algebras.

{\defn Let $(X,\mathcal O_X)$  be a $k$-scheme. The Lie-bracket defined  by the formula (\ref{Liebracket}) on each affine open subset defines \textit{a Lie-bracket} on the coherent sheaf $\Theta_{X/k} = \mathrm{Der}_k(\mathcal O_X)$ of derivations of the structural sheaf of $X$.

Similarly, if $(M,\mathcal O_M)$ is a (real or complex) analytic space, then the formula (\ref{Liebracket}) defines a Lie-bracket on the coherent sheaf $\Theta_{M} = \mathrm{Der}(\mathcal O_M)$. \\}
 
In both cases, under an additional smoothness assumption on $X$ (respectively on $M$), the sheaf $\Theta_{X/k}$ is a locally free on $X$ and is, indeed, the sheaf of sections of the vector bundle $T_{X/k}$ --- or in other words, the sheaf of vector fields on  $X$.

{\lem Let $k$ be either the field of real or complex numbers and $(X,\mathcal O_X)$  be a $k$-scheme. The (algebraic) Lie-bracket on $X$ satisfies the obvious compatibility relation with the analytic one on $X(k)^{an}$:

$$ [v,w]^{an} = [v^{an}, w^{an}]$$

where $-^{an}$ denotes either the real or the complex analytification and the Lie-bracket on the right-hand side is the (real or complex) analytic one.}

\begin{proof}
This can easily be derived from the standard properties of the analytification functor. It is also a direct consequence of the formula (\ref{coordinates-Lie-bracket}) of the next paragraph in both analytic and étale coordinates.    
\end{proof}

{\defn Let $(X,\delta_X)$ be a $D$-scheme over a constant differential field $(k,0)$ and $v$ a vector field on $X$. The \textit{Lie-derivative of $\delta_X$}, denoted $\mathcal L_{\delta_X} : \Theta_{X/k} \longrightarrow \Theta_{X/k} $ is the $k$-linear morphism defined by:
$$ \mathcal L_{\delta_X}(\delta) = [\delta_{X |U},\delta]$$
for every local section $\delta \in \Theta_{X/k}(U)$.}

{\lem\label{derivation}  Let $(X,\delta_X)$ be a $D$-scheme over a constant differential field $(k,0)$, let $\delta,\delta' \in \Theta_{X/k}(U)$ be two local sections defined on the same open set $U$ and $a \in \mathcal O_X(U)$. We have:
$$\begin{cases}
\mathcal L_{\delta_X}(\delta + \delta') = \mathcal L_{\delta_X}(\delta) + \mathcal L_{\delta_X}(\delta') \\
\mathcal L_{\delta_X}(a.\delta) = \delta_X(a). \delta + a. \mathcal L_{\delta_X}(\delta). 
\end{cases}$$}

\begin{proof}
These two properties follow immediately from the formula (\ref{Liebracket}).
\end{proof}

Before studying, more generally, the $k$-linear operator on a coherent sheaf satisfying the properties of Lemma \ref{derivation}, we compute, in the next paragraph, the Lie-bracket of two vector fields on a smooth variety in local coordinates (analytic coordinates in the analytic case and étale coordinates in the algebraic one).
  
\subsection{Computation in analytic and étale coordinates} Every analytic manifold $(M,\mathcal O_M)$ can be covered by analytic charts. More precisely, there exists a covering of $M$ by open subsets endowed with analytic coordinates $x_1,\ldots,x_n$ (meaning that the map $ x= (x_1,\ldots,x_n) : U \longrightarrow k^n$ is an analytic isomorphism onto its image). 

In order to prove local properties of the Lie-derivative, we will sometimes work locally inside these coordinates. Instead of analytic coordinates, we will use étale coordinates when working with algebraic varieties.
 
{\defn Let $U$ be a smooth variety over some field $k$ of dimension $n$. A \textit{system $(x_1,\ldots, x_n)$ of étale coordinates on $U$} is an étale morphism $x = (x_1,\ldots,x_n) : U \longrightarrow \mathbb{A}^n$.

In other words, it is a $n$-tuple $(x_1, \ldots , x_n)$ of regular functions on $U$ such that the section $dx_1 \wedge \cdots \wedge dx_n$ of the canonical line bundle $\Omega^n_{U/k} = \Lambda^n \Omega^1_{U/k}$ does not vanish on $U$. }

{\rem Note that, if $(x_1,\ldots x_n): U \longrightarrow \mathbb{A}^n$ is a system of étale coordinates then, by definition, $\lbrace dx_1, \ldots , dx_n\rbrace$  define a trivialization of $\Omega^n_{U/k}$. It follows that the dual basis of vector fields $\lbrace \frac \partial {\partial x_1},  \ldots , \frac \partial {\partial x_n} \rbrace$ define a trivialization of the tangent bundle of $X$.}

{\lem Let $X$ be a smooth algebraic variety over some field $k$ of characteristic $0$. There exists a covering of $X$ by Zariski-open subsets $(U;x_1,\ldots x_n)$ endowed with étale coordinates. \qed} 

{\exam Let $X$ be a smooth algebraic variety over a field $k$ of characteristic $0$. Consider an open set $U \subset X$ and $(x_1,\ldots x_n)$ a system of étale coordinates on $U$. Since the vector fields $\frac \partial {\partial x_1},  \ldots , \frac \partial {\partial x_n}$ define a trivialization of the tangent bundle of $X$, we may write

$$v_{|U} = \sum_{i = 1}^n v_i \frac \partial {\partial x_i} \text{ and } w_{|U} = \sum_{i = 1}^n w_i \frac \partial {\partial x_i}$$
for some functions $v_i,w_i \in \mathcal O_X(U)$.}

{\lem With the notation above, the Lie-bracket of the vector fields $v$ and $w$ is given by:
\begin{eqnarray}\label{coordinates-Lie-bracket}
 [v,w]_{|U} = \sum_{i = 1}^n \bigg ( \sum_{j = 1}^n v_j \frac {\partial w_i} {\partial x_j} - w_j \frac {\partial v_i} {\partial x_j} \bigg )  \frac \partial {\partial x_i}.
\end{eqnarray}}

\begin{proof}
The lemma follows from Lemma \ref{derivation} applied to both $\mathcal L_v$ and $\mathcal L_w$.
\end{proof}
Note that when one works with (a subfield of) the field of complex numbers, on may use analytic coordinates instead of \'etale coordinates and the formula (\ref{coordinates-Lie-bracket}) also holds in this analytic setting.

\subsection{Lie-derivative and $D$-coherent sheaves} The notion of $D$-coherent sheaves over some $D$-scheme $(X,\delta_X)$ formalizes the notion of a ``\textit{linear differential equation over} $(X,\delta_X)$''.                                                                                                                                                                                                                                                                                                                                                                                                                                                                                                                                                                                                                                                                                                                                                                                                                                                                                                                                                                                                                                                                                                                                                                                                                                                                                                                                                                                                                                                                                                                                                                                                                                                                                                                                                                                                                                                                                                                                                                                                                                                                                                                                                                                                                                                                                                                                                                                                                                                                                                                                                                                                                                                                                                                                                                                                                                                                                                                                                                                                                                                                                                                                                                                                                                                                                                                                                                                                                                                                                                                                                                                                                                                                                                                                                                                                                                                                                                                                                                                                                                         This is a straightforward generalization to $D$-schemes of the notion of $\delta$-module over differential fields, that appears for example in \citep[Section 3]{Pil}. 

{\defn Let $(X,\delta_X)$ be a $D$-scheme over some constant differential field $(k,0)$. A \textit{$D$-coherent sheaf} over $(X,\delta_X)$ is a pair $(\mathcal E, \nabla)$ where $\mathcal E$ is a coherent sheaf over $X$ and $\nabla : \mathcal E \longrightarrow \mathcal E$ is a $k$-linear sheaf morphism satisfying the Leibniz-rule with respect to scalar multiplication:
$$ \nabla(a.m) = \delta_X(a).m + a.\nabla(m)$$
for every local sections $a \in \mathcal O_X(U)$ and $m \in \mathcal E(U)$ on some open subset $U$ of $X$.

If $(\mathcal E,\nabla_\mathcal E)$ and $(\mathcal F,\nabla_\mathcal F)$ are both $D$-coherent sheaves over $(X,\delta_X)$, then \textit{a morphism of $D$-coherent sheaves over $(X,\delta_X)$} is a morphism $f : \mathcal E \longrightarrow \mathcal F$  of coherent sheaves over $X$ such that $$f \circ \nabla_E = \nabla_F \circ f.$$}

{\rem The notion of $D$-coherent sheaf is closely related to the more usual notion of a coherent sheaf $\mathcal E$ endowed with a connexion $\nabla$. Recall that if $X$ is a scheme over a field $k$, a \textit{connexion $\nabla$ on a coherent sheaf $\mathcal E$} over $X$ is $k$-bilinear morphism:
$$\nabla : \mathcal E \times \Theta_{X/k} \longrightarrow \mathcal E$$ }
which satisfies the Leibniz rule with respect to scalar multiplication on $\mathcal E$ and is $\mathcal O_X$-linear with respect to scalar multiplication on $\Theta_{X/k}$.

{\lem\label{connexion} Let $(X,\delta_X)$ be a $D$-scheme over some constant differential field $(k,0)$ and $(\mathcal E,\nabla)$ a coherent sheaf endowed with a connexion on $X$. Then $(\mathcal E, \nabla_{\delta_X})$ is a $D$-coherent sheaf.}

\begin{proof}
By definition, the morphism $\nabla_{\delta_X}$ is $k$-linear and satisfies the Leibniz rule.
\end{proof}
In particular, we get the following example: 

{\exam\label{trivialDstructure} Let $(X,\delta_X)$ be a $D$-scheme over some constant differential field $(k,0)$ and $\mathcal E = \mathcal O_X \epsilon_1 \oplus \cdots \oplus \mathcal O_X \epsilon_n$ be a free sheaf of rank $n$ over $X$. Define the $k$-linear map $\nabla_0 : \mathcal E \longrightarrow \mathcal E$ by the formula:
$$ \nabla_0 (\sum_{i = 1}^n f_i \epsilon_i) = \sum_{i = 1}^n  \delta_X(f_i) \epsilon_i.$$  
Then $(\mathcal E, \nabla_0)$ is a $D$-coherent sheaf over $(X,\delta_X)$.}

{\lem\label{affinestructure} Let $(X,\delta_X)$ be a $D$-scheme over some constant differential field $(k,0)$ and let $\mathcal E$ be a coherent sheaf on $X$. If $(\mathcal E,\nabla)$ and $(\mathcal E,\nabla')$ are both $D$-coherent sheaves then:
$$\nabla - \nabla' \in \mathrm{End_{\mathcal O_X}}(\mathcal E).$$} 

\begin{proof}
For a local function $a \in \mathcal O_X(U)$ and a local section $\sigma \in \mathcal O_X(U)$, we have:
$$ (\nabla - \nabla')(a.\sigma) = a.(\nabla(\sigma) - \nabla'(\sigma)) + \delta_X(a).\sigma - \delta_X(a).\sigma = a.(\nabla - \nabla')(\sigma).$$  
It follows that $\nabla - \nabla'$ is $\mathcal O_X$-linear.
\end{proof}
{\exam Let $(X,\delta_X)$ be a smooth $D$-variety over a constant differential field $(k,0)$ and $(\mathcal E,\nabla)$ a locally free $D$-coherent sheaf.

Consider an open set $U \subset X$ for which $\mathcal E_{|U} = \mathcal O_U \epsilon_1 \oplus \cdots \oplus \mathcal O_U \epsilon_n$ is free. Using Example \ref{trivialDstructure} and Lemma \ref{affinestructure}, there are functions $a_{i,j} \in \mathcal O_X(U)$ for $i,j = 1,\ldots, n$ such that:

\begin{equation}\label{normalform}
\nabla(\sum_{i = 1}^n f_i \epsilon_i) = \nabla_0(\sum_{i = 1}^n f_i \epsilon_i) + A.(f_1,\ldots, f_n)
\end{equation}
where $A = (a_{i,j})$ is the $n \times n$-matrix with coefficients $a_{i,j}(x) \in \mathcal O_X(U)$.

Conversely, for every matrix $A = (a_{i,j})_{i,j \leq n}$ with coefficients $a_{i,j}(x) \in \mathcal O_X(U)$, $\nabla = \nabla_0 + A$ defines a $D$-coherent sheaf on $\Theta_{U/k}$.  }

{\exam\label{LiederivativeasDcoherent} Let $(X,\delta_X)$ be a $D$-scheme over some constant differential field $(k,0)$. By Lemma \ref{derivation}, the pair $(\Theta_{X/k}, L_{\delta_X})$ is a $D$-coherent sheaf over $(X,\delta_X)$.

Assume now that $X$ is a smooth variety. The derivation $\delta_X$ is associated to a vector field $v$ on $X$. Consider $U \subset X$ an open subset and $(x_1,\ldots, x_n)$ a system of étale coordinates on $U$. Combining the formulas (\ref{normalform}) and (\ref{coordinates-Lie-bracket}), we get:
$$ L_v(w) = \nabla_0(w) - A.w \text{ where } a_{i,j} = \frac {\partial v_i} {\partial x_j} \in \mathcal O_X(U).$$
} 

Since the coefficients $a_{i,j}$ of the matrix $A$ do not depend linearly on the vector field $v$, this kind of $D$-coherent sheaves never come by Lemma \ref{connexion} from a connexion on $X$.  
 
\subsection{Categories of $D$-coherent sheaves over $D$-schemes} 
  
{\lem\label{abelian-category} Let $(X,\delta_X)$ be a $D$-scheme over some constant differential field $(k,0)$. The category of $D$-coherent sheaves over $(X,\delta_X)$  is an Abelian category.}

\begin{proof}
If $f : (\mathcal E,\nabla_\mathcal E) \longrightarrow (\mathcal F,\nabla_\mathcal F)$ is a morphism of $D$-coherent sheaves, it is easy to check that the kernel and the image of $f$ are respective $D$-subcoherent sheaves of $(\mathcal E,\nabla_\mathcal E)$ and $(\mathcal F,\nabla_\mathcal F)$ respectively.

Moreover, if $\mathcal G$ is a $D$-coherent sheave of $(\mathcal E,\nabla_\mathcal E)$, then there exists a unique $D$-coherent sheaf structure $(\mathcal E / \mathcal G, \nabla_{\mathcal E / \mathcal G})$ which makes the canonical projection into a morphism of $D$-coherent sheaves.

Using the forgetful functor to the category of coherent sheaves on $X$, which is an Abelian category, it is easy to check that the axioms of an Abelian category are satisfied. 
\end{proof}

{\defn Let $\phi : (X,\delta_X) \longrightarrow (Y,\delta_Y)$ be a morphism of $D$-schemes over some constant differential field $(k,0)$ and $(\mathcal E, \nabla_\mathcal E)$ a $D$-coherent sheaf over $(Y,\delta_Y)$.

The \textit{pull-back of $(\mathcal E, \nabla_\mathcal E)$ by $\phi$}, denoted  $\phi^\ast (\mathcal E, \nabla_\mathcal E)$, is the coherent sheaf $\phi^\ast \mathcal E$ over $X$ endowed with the derivation: 
$$ \phi^\ast \nabla_\mathcal E =  \nabla_\mathcal E \otimes 1 + \mathrm{Id} \otimes \delta_X.$$}

{\exam Let $\phi : (X,\delta_X) \longrightarrow (Y,\delta_Y)$ be a morphism of $D$-schemes over some constant differential field $(k,0)$ and $(\mathcal E, \nabla_\mathcal E)$ a locally free $D$-coherent sheaf over $(Y,\delta_Y)$.

Consider an open subset $U \subset Y$ such that the restriction $\mathcal E_{|U} = \mathcal O_U \epsilon_1 \oplus \cdots \oplus \mathcal O_U \epsilon_n$ of $\mathcal E$ to $U$ is free. Using formula (\ref{normalform}), we can write: 
$$ \nabla_{\mathcal E |U} = \nabla_0 + A$$  }
where $A = (a_{i,j})$ is an $n\times n$ matrix with coefficients $a_{i,j} \in \mathcal O_Y(U)$.

Set $V = \phi^{-1}(U)$. Note that the restriction $\phi^\ast \mathcal E_{|V}  = \mathcal O_V \phi^\ast \epsilon_1 \oplus \cdots \oplus \mathcal O_V \phi^\ast \epsilon_n$ of $\phi^\ast \mathcal E$ to $V$ is also free.

{\lem\label{basechange} With the notations above, the $D$-coherent structure on $\phi^\ast \mathcal E$ is given by: 
$$ \phi^\ast \nabla_\mathcal E = \nabla_0 + A^\phi$$
where $A^\phi$ is the matrix with coefficients $a_{i,j} \circ \phi \in \mathcal O_X(V).$}

\subsection{Main Proposition}

{\Prop\label{functoriality}  Let $\phi : (X,v) \longrightarrow (Y,w)$ be a morphism of smooth $D$-varieties over some constant differential field $(k,0)$. The derivative of $\phi$ defines a morphism of $D$-coherent sheaves over $X$:
$$d\phi :  (\Theta_{X/k},\mathcal L_v)  \longrightarrow \phi^\ast (\Theta_{Y/k} , \mathcal L_w)$$ }

\begin{proof}
It is sufficient to work locally in the Zariski topology. Consider \'etale coordinates $(x_1,\dots, x_p) : U \longrightarrow \mathbb{A}^p$ and $(y_1,\dots, y_q) : V \longrightarrow \mathbb{A}^q$ be \'etale coordinates on $Y$ such that $\phi : U \longrightarrow V$.

We want to show that $d\phi \circ \mathcal L_v = \phi^\ast \mathcal L_w \circ d\phi$ . Using Example \ref{LiederivativeasDcoherent} and Lemma \ref{basechange}, we can write: 
$$
\begin{cases}
\mathcal L_v = \nabla_0 - A \text{ where } a_{i,j} = \frac {\partial v_i} {\partial x_j} \in \mathcal O_X(U) \\
\phi^\ast \mathcal L_w = \nabla'_0 - B^\phi \text{ where } b_{i,j} = \frac {\partial w_i} {\partial y_j} \in \mathcal O_Y(V).
\end{cases}$$

Using this notations, the previous equality translates into:
$$ B^\phi . d\phi = d\phi. A   + (\nabla'_0. d\phi - d\phi.\nabla_0)$$
which is an identity between two matrices of size $p \times q$ with coefficients in $\mathcal O_X(U)$.

Now, since $\phi$ is a morphism of $D$-varieties, we have $d\phi(v) = \phi^\ast w$, which --- after denoting $\phi_j = y_j \circ \phi$, the coordinate function of $\phi$  --- translates in these coordinates by:
$$\sum_{k=1}^p \frac {\partial \phi_j}{\partial x_k} v_k = w_j \circ \phi.$$

For $1 \leq i \leq p$ and $1 \leq j \leq q$, the chain rule for derivation as well as the Leibniz rule imply that:
$$ (B^\phi . d\phi)_{i,j} = \frac \partial {\partial x_i} (w_j \circ \phi) = \frac \partial {\partial x_i} (\sum_{k = 1}^p  \frac {\partial \phi_j}{\partial x_k} v_k) = (d\phi.A)_{i,j} + \sum_{k = 1}^p v_k \frac {\partial^2 \phi_j}{\partial x_k x_i}.$$

Moreover, since $\nabla_0 (\frac \partial {\partial x_i}) = 0$, we have that:

$$ (\nabla'_0. d\phi - d\phi.\nabla_0)_{i,j} = (\nabla'_0. d\phi)_{i,j} = v( \frac {\partial \phi_j} {\partial x_i}) =  \sum_{k = 1}^p v_k \frac {\partial^2 \phi_j}{\partial x_k x_i}.$$  
This concludes the proof of the proposition.
\end{proof}

We now gather two corollaries of Proposition \ref{functoriality} which deal respectively with two different geometric situations.
 
{\cor\label{corollary1} Let $\phi : (X,v) \longrightarrow (Y,w)$ be a dominant morphism of smooth $D$-varieties over a constant differential field $(k,0)$. The coherent subsheaf $\Theta_{X/Y} = \mathrm{Ker}(d\phi)$ is invariant under the Lie-derivative $L_v$ of $v$.}

{\cor Let $(X,v)$ be a smooth $D$-variety and $Y \subset_i X$ a closed smooth invariant submanifold. We have an exact sequence of sheaves over $Y$:
$$ 0 \longrightarrow \Theta_{Y/k} \longrightarrow i^\ast \Theta_{X/k} \longrightarrow \mathcal N_{X/Y} \longrightarrow 0$$
where $\mathcal N_{X/Y}$ denotes the normal bundle of $Y$ in $X$. The Lie-derivative $\mathcal L_v$ induces a well-defined $D$-coherent sheaf structure on $\mathcal N_{X/Y}$.}

\subsection{Cauchy formula on an analytic manifold}

{\cons Let $(A,\delta_A)$ be a differential ring of characteristic $0$. Consider the differential ring $(A((t)), \frac d {dt})$ of formal power series over $A$ endowed with the derivation sending $t$ to $0$ and the morphism of differential rings $ \phi^\ast :  (A,\delta_A) \longrightarrow (A((t)), \frac d {dt})$  given by the expansion in power series:
$$ a \mapsto \sum_{k=0}^\infty \frac {\delta_A^k(a)} {k!} t^k.$$

We denote by $-_{|t = 0}:  A((t)) \longrightarrow A$ the morphism of evaluation at $t = 0$ and we denote by $\delta_A$, the unique extension of $\delta_A$ to $A((t))$ satisfying $\delta_A(t) = 0$. Note that with this derivation, the morphism of rings 
 $-_{|t = 0}$ becomes a morphism of differential rings:
 
{\lem\label{lemmapowerseries} With the notation above for every $f \in A$, $\phi^\ast(f)$ may be described as the unique solution of the differential equation: 
$$ \begin{cases}
\phi^\ast(f)_{|t = 0} = f \\
\frac d {d t} \phi^\ast(f) = \delta_A (\phi^\ast(f)) 
\end{cases} $$}}

{\rem Consider $M$ an analytic manifold, $v$ an analytic vector field on $M$ and $a \in M$. The vector field $v$ induces a derivation $\delta_v$ on $A = \mathcal O_{M,a}$. In \citep[Lemme 3.1.21]{moi}, we proved --- using Cauchy integral formula the bound the norm of $\delta_v^n(f)$ --- that in that case, the morphism $\phi^\ast$ may be factored as:  
$$ \phi^\ast :  A \longrightarrow \mathcal O_{M \times \mathbb{C}, (a,0)} \subset A((t)).$$

As noted in \cite{moi}, one easily checks that if $\phi_a$ denotes the local analytic flow of the vector field at $v$, this implies that: 
$$ f \circ \phi_a = \phi^\ast(f) = \sum_{k = 0}^\infty \frac {\delta_v^k(f)} {k!} t^k. $$
This property was then used to translate --- for a closed submanifold of $M$ --- invariance properties with respect to the vector field in terms of invariance properties with respect to the local analytic flow (see \cite[Proposition 3.1.20]{moi}).}

\cons Let $(A,\delta_A)$ be a differential ring of characteristic $0$ and $(M,\nabla_M)$ be a $D$-module over $A$. One may define a morphism of $A$-modules $M \longrightarrow M \otimes A((t))$ by the formula:
$$ m \mapsto \sum_{k = 0}^\infty \frac {\nabla_M^k(m)} {k!} \otimes t^k.$$ 

When  $(M,\nabla_M)$ is given by the Lie-derivative of a vector field $v$ on an analytic manifold $M$, we will use once again Cauchy integral formula to bound the sequence $(\frac {\nabla_M^k(m)} {k!})_{k \in \mathbb{N}}$:

Let $M$ be a smooth analytic manifold and let $v$ be a vector field on $M$.  For every point $a \in M$, the Lie-derivative of vector field $v$ defines a derivation $L_v : \Theta_{M,a} \longrightarrow \Theta_{M,a}$. We denote by $\mathcal O_{M,a}\lbrace t \rbrace$ the ring of local analytic functions on $\mathbb{C} \times M$. Pull-back along the local flow $\phi_a$ of $v$ defines a morphism of $\mathcal O_{M,a}$-modules: 
$$ \phi_a^\ast : \Theta_{a,M} \longrightarrow \Theta_{a,M} \otimes \mathcal O_{M,a}\lbrace t \rbrace$$

We extend the derivation $\mathcal L_v$ on $\Theta_{a,M}$ to a derivation on  $\Theta_{a,M} \otimes \mathcal O_{M,a}\lbrace t \rbrace$  still denoted $\mathcal L_v$ by setting $\mathcal L_v(t) = 0$.

{\lem[Cauchy formula for the Lie derivative]\label{Cauchyformulas} With the notation above for every vector fields $w$ in a neighborhood of $p$, we have:
$$\begin{cases}
\phi_a^\ast(w)_{|t = 0} = w \\
\frac d {dt} \phi_a^\ast(w) = \mathcal L_v(\phi_a^\ast(w))
\end{cases}$$
Moreover, we have: 
$$  \phi_a^\ast w = \sum_{n = 0}^\infty \frac {\mathcal L_v^n(w)} {n!} t^n.$$}

\begin{proof}
The relation between the Lie-derivative of a vector field $v$ and the pull-back by the local flow of the vector field $v$, given in the first part is well-known and holds more generally on smooth manifolds (see any textbook of differential geometry).

For the second part, one needs to check that the right hand-side converges normally to an analytic function in a neighborhood of $p$. By formal derivation of power series, it follows easily that the left-hand side satisfies the differential equation given by the first part and therefore must be equal to the local flow $\phi_a^\ast w$. Similar formulas already appear in the work of Cauchy.

Given a vector field $w$ on $M$, we need a uniform bound for $\mathcal L^n_v(w)$ in a neighborhood of $p$, in order to prove normal convergence. Since $M$ is smooth and that we work locally on $M$, we may assume that $M = \mathbb{C}^n$ and $p = 0$. 

Fix $r_0 < R_0$ are the radii of two complex polydisks of $p$ and $n$ a natural number. 
By applying Cauchy integral formula to the holomorphic coordinates of $w$, there exists a constant $C > 0$ such that for all $\epsilon > 0$, all radii $r > 0$ and all vector fields $w$, we have:  
$$ || \mathcal L_v(w) ||_{\infty,r - \epsilon} \leq  \frac C \epsilon ||w||_{\infty, r}.$$ 

where $||-||_{\infty, r}$ denotes the supremum norm on the polydisk with radius $r$ . It then follows from $n$ successive applications of the previous inequality to the radii $r_k = r_0 + (R_0 - r_0)\frac k n$  that:

$$ || \mathcal L^n_v(w) ||_{\infty,r_0}  \leq  (\frac {C} {n})^n ||w||_{\infty, R_0} $$

The normal convergence of the left-hand side follows from this inequality.
\end{proof}

\section{Foliations on a smooth algebraic variety}

In this section, we recall the standard definition of an algebraic foliation --- in the setting of algebraic varieties over a field $k$ of characteristic $0$ --- that we will use in this article.

Let $X$ be a smooth irreducible algebraic variety over $k$. Intuitively, a foliation $\mathcal F$ on $X$ is the data of a subspace $T_{\mathcal F,x}$ of the tangent space $T_{X,x}$ at $x$ for every point $x$ of $X$ that depends algebraically on the point $x \in X$ and such that the sheaf of sections of $\mathcal F$ is stable under Lie-bracket.

The dimension of the subspace $F_\eta$ at the generic point $\eta$ of $X$ is called the \textit{rank of the foliation}. A \textit{singularity of the foliation} $\mathcal F$ is simply a point $x$ of $X$ where the dimension of the fibre  $F_x$ is less than the rank of the foliation.

More restrictively, we will require the foliations to satisfy an additional assumption of \textit{saturation}. On the one hand, the involutivity property --- that is the stability under Lie-bracket --- ensures the local analytic integrability of the foliation, \textit{outside of the singular locus}. On the other hand, the saturation hypothesis ensures that the singular locus of $\mathcal F$ is \textit{``small''}, namely of codimension at least $2$ in $X$.

The main motivation for these additional requirement is the extension result (Proposition \ref{saturation}) for any algebraic foliation $\mathcal F$ on a quasi-projective smooth variety $X$ to a foliation $\overline{\mathcal F}$ on any projective closure.

\subsection{Algebraic foliations and their singular locus} Let $X$ be an irreducible and smooth variety over a field $k$ of characteristic $0$. The coherent sheaf $\Theta_{X/k}$ is a locally free sheaf of rank $n = \mathrm{dim}(X)$.

{\defn\label{foliation} A \textit{foliation} $\mathcal F$ on $X$ is a sub-sheaf $T_\mathcal F$ of the tangent bundle $\Theta_{X/k}$ which satisfies:
\begin{itemize}
\item[(i)] The sub-sheaf $T_\mathcal F$ is \textit{involutive}, that is, stable under the Lie-bracket.
\item[(ii)] The sub-sheaf $T_\mathcal F$ is \textit{saturated}, that is, the quotient $\Theta_{X/k} / T_\mathcal F$ does not have torsion.
\end{itemize}
The coherent sheaf $T_\mathcal F$ is \textit{the sheaf of vector fields tangent to the foliation $\mathcal F$}. The \textit{rank of the foliation $\mathcal F$} is the (generic) rank of the coherent sheaf $T_\mathcal F$. Alternatively, we say that $\mathcal F$ is a $k$-foliation to mean that $\mathcal F$ is a foliation with rank $k$.}

{\rem We first comment on the two main assumptions in Definition \ref{foliation}, namely involutivity and saturation.

\begin{itemize}
\item[(1)] In this article, we are mainly interested in smooth algebraic varieties defined over the field of real or complex numbers. In that case, the \textit{involutivity} assumption is crucial to ensure local analytic integrability of the algebraic foliation under study (see section 2.5). However, for a coherent sub-sheaf $\mathcal L$ of $\Theta_{X/k}$ of rank $1$, the condition of involutivity is automatically satisfied.

\item[(2)] The property of \textit{saturation} for $T_\mathcal F$ (inside a locally free sheaf) implies that the coherent sheaf $T_\mathcal F$ is reflexive (namely, isomorphic to its bidual). This property, which is extensively studied in \cite{Har2}, is weaker than being locally free but stronger than being torsion-free.

In particular, a $1$-foliation is always defined by an invertible sheaf $T_\mathcal F$ whereas for a $2$-foliation on a smooth algebraic variety $X$ of dimension $3$  is always defined by a coherent sheaf $T_\mathcal F$ which is locally free outside a finite set of points (see \cite{Har2}) of $X$. Since the main result of this article deals with $D$-variety of dimension $3$, these two examples are the most important ones for the proof of Theorem \ref{maintheoremInt}.

\end{itemize}}

{\lem Let $X$ be a smooth algebraic variety over $k$ and let $T_\mathcal F \subset \Theta_{X/k}$ be a coherent sub-sheaf of the tangent sheaf of $X$. The following properties are equivalent: 
\begin{itemize}
\item[(i)] $\Theta_{X/k} / T_\mathcal F$ does not have torsion.
\item[(ii)] There exists an open set $U \subset X$ such that $codim_X(X \setminus U) \geq 2$ and 
$$ 0 \longrightarrow T_{\mathcal F|U} \longrightarrow \Theta_{U/k} \longrightarrow \Theta_{U/k} / T_{\mathcal F |U} \longrightarrow 0$$
is an exact sequence of locally free-sheaves on $U$.
\end{itemize}}

\begin{proof}
Indeed, fix $x \in X$ of codimension $1$. Since $X$ is a smooth algebraic variety, the local ring $\mathcal O_{X,x}$
is a principal local ring. It follows that the exact sequence of torsion-free $\mathcal O_{X,x}$-modules:
$$ 0 \longrightarrow T_{\mathcal F,x} \longrightarrow \Theta_{X/k,x} \longrightarrow \Theta_{X/k,x} / T_{\mathcal F,x} \longrightarrow 0$$ 
is in fact an exact sequence of free $\mathcal O_{X,x}$-module which therefore have to split. Since this is true for every point   $x \in X$ of codimension $1$, the lemma follows.
\end{proof}

{\defn\label{singularlocus} Let $\mathcal F$ be a foliation  on $X$. The \textit{singular locus of $\mathcal F$}, denoted $\mathrm{Sing}(\mathcal F)$ is the  set of points $x \in X$ such that $\Theta_{X/k}/\mathcal F$ is not locally-free in a neighborhood of $x$, i.e.
$$\mathrm{Sing}(\mathcal F) = \lbrace x \in X \text{ | } (\Theta_{X/k}/\mathcal F)_x \text{ is not a  }\mathcal O_{X,x}-\text{free module} \rbrace.$$}

{\defn\label{nonsingular} Let $\mathcal F$ be a foliation on $X$. The foliation $\mathcal F$ on $X$ is called \textit{non-singular} if $\mathrm{Sing}(\mathcal F) = \emptyset$.}

{\Prop Let $\mathcal F$ be a foliation  on $X$. The singular locus of $\mathcal F$ is a closed subset of $X$ of codimension $\geq 2$.}

\begin{proof}
The fact that the singular locus of a foliation is closed follows from the general properties of morphisms of coherent sheaves. This also follows from the computations of Example \ref{singularlocus-computation}. Moreover, the previous lemma shows that its codimension is greater than $2$.
\end{proof}

{\exam\label{singularlocus-computation} Let $\mathcal F$ be a foliation  on $X$ and $U \subset X$ an open subset endowed with a system of étale coordinates $x_1,\ldots x_n : U \longrightarrow \mathbb{A}^n$. Recall that the tangent sheaf $\Theta_{U/k}$ may be identified with the free sheaf $\mathcal O_{U/k}^n$ with basis $(\frac \partial {\partial x_1}, \ldots ,\frac \partial {\partial x_n})$.

Since the coherent sub-sheaf $\mathcal F_{|U}$ is finitely generated, it can be described as the coherent sheaf on $U$ spanned by some vector fields $v_1,\ldots, v_r$ written as:

$$ v_i = \sum_{i = 1}^n f_{i,j} \frac \partial {\partial x_i}.$$

We can now explicit the system of equations describing the singular locus of $\mathcal F$ restricted to $U$ in terms of these datas:
$$\mathrm{Sing}(\mathcal F) \cap U = \lbrace p \in U \text{ | } \mathrm{rank}(f_{i,j}(p)) \text{ is not maximal} \rbrace.$$  

Hence, if $\mathcal F$ is a foliation of rank $p$, the system of equations describing the singular locus of $\mathcal F$ is given by the vanishing of the $p \times p$ minors of the matrix $(f_{i,j})_{1 \leq i \leq r, 1 \leq j \leq n}$.}

\subsection{Analytification of an algebraic foliation} Assume that $k$ is the field of real or complex number. We have the following counterpart for Definition \ref{foliation} in the analytic setting:

{\defn\label{analyticfoliation} Let $M$ is a real or complex analytic manifold. An \textit{analytic foliation} on $M$ is a coherent subsheaf $\mathcal E$ of the tangent sheaf $\Theta_M$ of $M$ which is both involutive and saturated.}  

Similarly to the algebraic case, for an analytic foliation $\mathcal{E}$ on an analytic manifold $M$, one can define its \textit{rank} and its \textit{singular locus} $\mathrm{Sing}(\mathcal E)$. The singular locus of $\mathcal{E}$ is a closed analytic subspace of $M$ of codimension $\geq 2$.

{\lem Let $X$ be a complex (resp. a real) smooth algebraic variety and $\mathcal F$ a foliation on $X$. Through the complex-analytication (resp. real-analytification) functor $-^{an}$, $\mathcal F$ defines an analytic foliation $\mathcal F^{an}$ on $X^{an}$.

Moreover, the rank of $\mathcal F^{an}$ is the rank of $\mathcal F$ and the singular locus of $\mathcal F^{an}$ is the analytification of the singular locus of $\mathcal F$.}  

{\rem When we will work in the analytic setting, we will mainly work with \textit{non-singular foliation}. The only exception to that rule is the proof of Proposition \ref{invariance-singular} where we work with the local analytic flow of a vector field to prove invariance properties for its singular locus.

Once this has been established, we will simply throw away the singular locus $\mathrm{Sing}(\mathcal F)$ of the foliation $\mathcal F$ at the level of scheme before applying the analytification functor.}

\subsection{Saturation of an algebraic foliation on a open subset}

{\Prop[Saturation]\label{saturation} Let $X$ be a smooth algebraic variety over $k$ and let $U$ be a dense open subset.
Any algebraic foliation $\mathcal F$ on $U$ extends uniquely to an algebraic foliation on $X$.}

\begin{proof}

Let  $\mathcal F$ be an algebraic foliation on $U$. There exists a coherent subsheaf $\mathcal G$ of $\Theta_{X/k}$ such that $\mathcal G_{|U} = \mathcal F$ (namely, the sheaf $\mathcal G$  of vector fields $v$ such that $v_{|U} \in \mathcal F$, which, as a simple verification shows, is quasi-coherent, hence coherent).

Let $\overline{\mathcal F}$ be the saturation of $\mathcal G$ in $\Theta_{X/k}$ (see \cite{Har2}). By definition, $\overline{\mathcal F}$ is a saturated subsheaf of $\Theta_{X/k}$, whose restriction to $U$ is $\mathcal F$. Therefore, it suffices to check that $\overline{\mathcal F}$ is involutive.

Let $v\in \overline{\mathcal F}(V)$ be a local section of $\overline{\mathcal F}$. The Lie-bracket with $v$ defines a morphism of $\mathcal O_X$-modules:
$$ [ - , v ] : \overline{\mathcal F}_{|V} \longrightarrow (\Theta_{X/k} /\overline{\mathcal F})_{|V}$$

Since the algebraic foliation $\mathcal F$ is involutive, this morphism is zero on $U \cap V$. Consequently, the image of this morphism is a coherent sheaf whose support is a proper closed subvariety of $X$, hence a torsion sheaf. Since the coherent sheaf $\Theta_{X/k} /\overline{\mathcal F}$ has no torsion, this morphism is zero.

The uniqueness is a direct consequence of the saturation hypothesis, since two saturated subsheaves of $\Theta_{X/k}$ which have the same generic fibre are equal (see \cite{Har2}). 
\end{proof}
Proposition \ref{saturation} is useful to construct algebraic foliations on smooth algebraic variety $X$: it shows that one only needs to construct them on a dense open subset.

{\cons \label{rationalfactor} Let $\phi: X \dashrightarrow Y$ be a rational dominant morphism of smooth irreducible varieties over a field $k$. Denote by $n$ the dimension of $X$ and $m$ the dimension of $Y$.
Let $U$ be the biggest open set where $\phi$ is defined and smooth. The restriction of the differential $d\phi : \Theta_{X/k} \longrightarrow \phi^\ast \Theta_{Y/k}$ of $\phi$ to $U$ has constant rank $n-m$.

{\lem With the notation above, $\mathrm{Ker}(d\phi_{|U})$ defines a non-singular algebraic foliation on $U$ of rank  $(n-m)$.}

\begin{proof}
Since $\mathrm{Ker}(d\phi_{|U})$ is a sub-vector bundle of $T_{U/k}$, it suffices to prove that $\mathrm{Ker}(d\phi_{|U})$ is involutive.
Let $v$ be a local section of $\mathrm{Ker}(d\phi_{|U})$ on a open set $V$. Then, $\phi_{|V} : (V,v) \longrightarrow (Y,0)$ is a morphism of smooth $D$-varieties. By Corollary \ref{corollary1}, $\mathrm{Ker}(d\phi_{|U})$ is stable under $\mathcal L_v$.
Since this is true for every local section, $\mathrm{Ker}(d\phi_{|U})$ is involutive. 
\end{proof}

{\defn Let $\phi: X \dashrightarrow Y$ be a rational dominant morphism of smooth irreducible varieties over a field $k$. Denote by $n$ the dimension of $X$ and $m$ the dimension of $Y$.

By Proposition \ref{saturation}, the non-singular foliation  $\mathrm{Ker}(d\phi_{|U})$ uniquely extends to a (possibly singular) $(n-m)$-foliation on $X$ denoted $\mathcal F_\phi$ and called the \textit{foliation tangent to $\phi$}.}

\subsection{$1$-Foliation tangent to a vector field}

{\defn\label{rationalvectorfield} Let $X$ be a smooth and irreducible algebraic variety  and $v$ be a non-zero rational vector field on $X$.
Let $U$ be the biggest open set of $X$ where $v$ is defined and does not vanish. The restriction of $v$ on $U$ generates a non-singular $1$-foliation on $U$. By Proposition \ref{saturation}, this foliation uniquely extends to a $1$-foliation on $X$ denoted $\mathcal F_v$ and called the \textit{foliation tangent to $v$} (with possible singularities outside of $U$).}

{\lem Let $v,w$ be two non-zero rational vector fields on $X$. The foliations $\mathcal F_v$ and $\mathcal F_w$ coïncide if and only if $v = f.w$ for some non-zero rational function $f$ on $X$.}

\begin{proof}
The main observation is that the vector fields $v$ and $w$ define the same foliation if and only if they define the same foliation on an open set of $X$ (by uniqueness of Proposition \ref{saturation}).

Now, if $f$ is a non-zero rational function such that $v = f.w$ then the foliation tangent to $v$ and $w$ agrees at least outside of the indeterminacy locus and the zero locus of $f$. By the previous observation, they agree everywhere. Conversely, if $v$ and $w$ define the same foliation then, outside the closed subset of $X$ given as the union of the singularities and zeros of both $v$ and $w$, they differ by multiplication by a non vanishing function. Hence they differ by multiplication by a non-zero rational function.   
\end{proof} 
 
{\rem The main content of Proposition \ref{saturation} is that --- in contrast with vector fields --- foliations extend automatically to any smooth projective model of a quasi-projective variety.

If $(X,v)$ is any quasi-projective $D$-variety, one can consider the $1$-foliation $\mathcal F_v$ tangent to $v$ on the projective closure $\overline{X}$ of $X$. In particular, in that setting, the notion of \textit{closed invariant subvarieties} naturally \textit{extends} to closed subvarieties of $\overline{X}$.}

 We illustrate the previous remark by computing the foliation on $\mathbb{P}^2$ associated to a polynomial vector field on $\mathbb{A}^2$:
 
{\exam Let $k$ be a field of characteristic $0$. We consider polynomial vector fields $v$ on the plane $\mathbb{A}^2$:
$$v(x,y) = a(x,y)\frac \partial {\partial x} + b(x,y)\frac \partial {\partial y}$$
where $a(x,y)$ and $b(x,y)$ are two coprime polynomials with coefficients in $k$, with the same degree $n$. The condition of coprimeness simply express that $v $ defines a saturated subsheaf of $\mathbb{A}^2_{x,y}$.

We now describe the foliation $\mathcal F_v$ tangent to $v$ on $\mathbb{P}^2$ in the chart $(s,t)$ where:
$$s = \frac 1 x \text{ and } t = \frac y x.$$

Note that in the coordinates $(s,t)$, the line at infinity is described by $s = 0$. Denote by $\hat a(s,t) = s^n. a(1/s,t/s) \in k[s,t]$ and $\hat b(s,t) = s^n.b(1/s,t/s) \in k[s,t]$. A simple computation shows the coordinates $(s,t)$, the vector field $v$ is given by: 
$$s^{n-1}. v(s,t) = -s. \hat a(s,t)\frac \partial {\partial s} + (-t. \hat a(s,t)  + \hat b(s,t)) \frac \partial {\partial t}.$$

Now, if we denote by $w(s,t)$ the vector field on the right-hand side, the vector fields $w(s,t)$ and $v(s,t)$ have the same singularities outside of the line $x = 0$ and $s = 0$. Hence, the vector field $w(s,t)$ defines a saturated sub-sheaf of $\Theta_{\mathbb{A}^2_{s,t}/k}$ if and only if the vector field $w(s,t)$ does not vanish on the line at infinity. This condition can be rewritten as: 

$$ P(t) := -t. \hat a(0,t)  + \hat b(0,t) \neq 0 \text{ or }  Q(x,y) = x^nP(y/x) =  x.b_n(x,y) - y.a_n(x,y) \neq 0.$$ 

where $a_n(x,y)$ and $b_n(x,y)$ are respectively the homogeneous parts of degree $n$ of $a(x,y)$ and $b(x,y)$ respectively.

We have proven the following:

{\lem Let $v(x,y) = a(x,y)\frac \partial {\partial x} + b(x,y)\frac \partial {\partial y}$ be a polynomial vector field on the plane where $a(x,y)$ and $b(x,y)$ are two coprime polynomials with degree $n$. Assume moreover that: 
$$ Q(x,y) = x.b_n(x,y) - y.a_n(x,y) \neq 0$$
Then the line at infinity is invariant for the foliation $\mathcal F_v$ and the singularities $\mathrm{Sing}(\mathcal F_v)_\infty$ of $\mathcal F_v$ on the line at infinity are given (in homogeneous coordinates) as the roots of $Q(x,y)$.}

{\rem This simple lemma admits strong consequences for the structure of closed invariant varieties for the $D$-variety $(\mathbb{A}^2,v)$ for a vector field $v$ satisfying the assumption of the lemma.
For example, if $L_\infty$ denotes the line at infinity then every closed invariant curve $C$ of $(\mathbb{A}^2,v)$ satisfies: 

$$ C_\infty = \overline{C} \cap L_\infty \subset \mathrm{Sing}(\mathcal F_v)_\infty.$$

This property heavily constrains the possible closed invariant curves for $(\mathbb{A}^2,v)$. For example, with refinements of this idea in \cite{Coutalg}, Coutinho and Menasch\'e describe an algorithm to compute polynomial vectors fields on the plane without closed invariant curves.}

\subsection{Analytic leaves of a non-singular analytic foliation}

{\Thm[Frobenius Integrability Theorem, {\citep[Theorem 2.9, Part I]{Ily}}]\label{Frobeniustheorem} Let $M$ be a complex analytic manifold of dimension $n$, $\mathcal F$ be a non-singular analytic $r$-foliation on $M$ and $x \in M$.

There exist a neighborhood $U$ of $x$ in $M$, open subsets $V_1 \subset \mathbb{C}^r$, $V_2 \subset \mathbb{C}^{n-r}$ and an analytic isomorphism $\phi: U \longrightarrow V_1 \times V_2$ such that for all $t \in V_2$, the tangent space of $L_t = \phi^{-1}(V_1 \times \lbrace t \rbrace)$ at every point $a \in L_t$ is $F_a$.}

\rem Let $M$ be a complex analytic manifold of dimension $n$, $\mathcal F$ be a non-singular analytic $r$-foliation on $M$ and $x \in M$.

The germ at $x$ of the analytic subvariety $L_y = \phi^{-1}(V_1 \times \lbrace y \rbrace)$ where $y = \phi_2(x)$ does not depend on the chosen trivialisation $\phi$ given by Theorem \ref{Frobeniustheorem} and is called \textit{the germ of the leaf of $\mathcal F$ through $x$}.
 
{\defn Let $M$ be a complex analytic manifold of dimension $n$ and $\mathcal F$ be a non-singular analytic $r$-foliation on $M$. 
A \textit{leaf of $\mathcal F$} is a maximal immersed analytic subvariety $N$ of $M$, which is the germ of a leaf at any point of $N$.}

{\cor  Let $M$ be a complex analytic manifold of dimension $n$ and $\mathcal F$ be a non-singular analytic $r$-foliation on $M$. The leaves of $\mathcal F$ define a partition of $M$ into immersed analytic subvariety of dimension the rank of $\mathcal F$.}
\begin{proof}
Indeed, by maximality, two distinct leaves of $\mathcal F$ never intersect. Moreover, by Theorem \ref{Frobeniustheorem}, any point $x \in M$ is contained in a leaf.
\end{proof}

{\defn Let $X$ be an algebraic manifold and $\mathcal F$ a possibly singular foliation of $X$ . Denote by $U$ the complementary of the singular locus $\mathrm{Sing}(F)$ of $\mathcal F$.

A leaf $\mathcal L$ of the non-singular foliation $\mathcal F^{an}_{|U}$ is called \textit{algebraic} if there exists a closed algebraic subvariety $Z$ of $X$ such that 
$$Z(\C) \cap U(\C) = \mathcal L.$$}

{\exam\label{example} Let $X$ be a smooth complex algebraic variety. By \cite{Per}, there exists at least one rational vector field $v$ (in fact, most of them), such that the foliation $\mathcal F_v$ defined in Example \ref{rationalvectorfield} does not admit any algebraic leaf.

However, by definition, the foliation associated to a rational factor $\mathcal F_\phi$ in Example \ref{rationalfactor}, there exists a non-empty open set $U$ such that any leaf that encounters $U$ is algebraic.}

\subsection{Continuous foliations on an analytic manifold} Let $M$ be a (real or complex) analytic manifold and let $E \subset TM$ be \textit{a continuous distribution on $M$}, that is a continuous sub-bundle $E$ of $TM$. Under such a weak regularity assumption,  the property of \textit{involutivity} formulated in terms of the Lie-derivative is not available anymore to distinguish foliations (satisfying local-integrability properties)  from other distributions.

Let's start by recalling the definition of an Anosov flow, which evidences two important continuous distributions.
  
{\defn Let $(M,g)$ be a Riemannian manifold of dimension $\geq 3$ and $v$ a $\mathcal C^\infty$-vector field on $M$ with a complete flow $(\phi_t)_{t \in \R}$. The flow $(M , (\phi_t)_{t \in \R})$ is called \textit{an Anosov flow} if the vector field $v$ does not vanish and there exists a splitting of the tangent bundle $TM$ into continuous sub-bundles
\begin{eqnarray}\label{Anosovsplitting}
TM = E^{ss} \oplus \R.v \oplus E^{su} \end{eqnarray}
satisfying:
\begin{itemize}
\item[(i)] The sub-bundles $E^{ss}$ and $E^{su}$ are non-trivial bundles which are $(d\phi_t)_{t \in \R}$-invariant.
\item[(ii)] There exist $C,C' > 0$ and $0 < \lambda < 1$ such that for all $u \in E^{ss}$,
$$|| d\phi_t(u) || \leq  C. \lambda^t ||u|| \text{ and } || d\phi_{-t}(u) || \geq  C'. \lambda^{-t} ||u|| \text{ for all t > 0}.$$
\item[(iii)] There exist $C,C' > 0$ and $0 < \lambda < 1$ such that for all $w \in E^{su}$,
$$|| d\phi_t(w) || \geq  C. \lambda^{-t} ||w|| \text{ and } || d\phi_{-t}(w) || \leq  C'. \lambda^{t} ||w|| \text{ for all t > 0}.$$
where the norm on the tangent bundle $TM$ is given by the Riemannian metric on $M$. 
\end{itemize}}

The distributions $E^{ss}$ and $E^{su}$ are called respectively \textit{strongly stable and strongly unstable distributions}.
Note that if $M$ is compact then any two Riemannian metrics on $M$ are equivalent. Hence, for a vector $v$ on a smooth compact manifold $M$, the property that expresses that its real-analytic flow is Anosov is independent of the choice of a Riemannian metric on $M$. In particular, for compact manifold $M$, these two distributions are uniquely determined by the vector field $v$.

{\rem  The work of Hasselblat  in \cite{Hasselblatt} (that goes back to Anosov in the case of diffeomorphisms) imply that these distributions do not, in general, satisfy stronger regularity properties (such as being $\mathcal C^2$). Consequently, when working with Anosov flows, we will need to work at this degree of regularity.

The procedure that allows that --- which produces, after real analytication, from an algebraic foliation $\mathcal F$, a continuous foliation \textit{outside of the singular locus} --- will be essential for our purposes.}

{\defn\label{continuous} Let $M$ be a analytic manifold of dimension $n$.  A \textit{continuous (non-singular) foliation $\mathcal F$ on $M$ of rank $p$} is a continuous atlas $\lbrace (U_i, \phi_i) \text{ , } i \in I \rbrace$ on $M$ such that: 
\begin{itemize}
\item[(i)] The image $\phi_i(U_i) = V_i^1 \times V_i^2 \subset \mathbb{R}^{p} \times \R^{n-p}$ can be written as a product of two connected open subsets $V_i^1$ and $V_i^2$ of $\R^p$ and $\R^{n-p}$ respectively. 
\item[(ii)] The transition maps $\psi_{i,j}$ are ``locally triangular'' with respect to the decomposition $\R^n = \R^p \oplus \R^{n-p}$ i.e.
$$ \psi_{i,j}(x,y) = (h_{i,j}(x,y), g_{i,j}(y)) \text{ where } x \in \R^p \text{ and } y \in \R^{n-p}.$$ 
for some functions $h_{i,j}$ and $g_{i,j}$ defined on the appropriate open subspaces of $\R^n$.
\end{itemize}

The \textit{leaves of the foliation $\mathcal F$} are then the equivalence classes for the equivalence relation $R$ generated by $x R y$ if: 
\begin{itemize}
\item[(iii)] $\exists i \in I$, $x,y \in U_i$  and the last $(n-p)$  coordinates of $\psi_i(y)$ and $\psi_i(x)$ are equal.
\end{itemize}}

In particular, the leaves $\lbrace L_\alpha , \alpha \in A \rbrace$ of the foliation $\mathcal F$ define a partition of $M$. The compatibility between this definition and Definition \ref{analyticfoliation} in the analytic case is ensured by Frobenius Integrability Theorem.

{\Thm[Hadamard-Perron Theorem -- {\citep[Section 2.2.i]{Hass}}] Let $(M,(\phi_t)_{t \in \R})$ be an Anosov flow as above. The stable and unstable distributions $E^s$ and $E^u$ are integrable by continuous foliations with $\mathcal C^1$ leaves.}

\section{Foliations invariant by a vector field}

In this section, we study foliations in the setting of smooth $D$-variety $(X,v)$ over some constant differential field $(k,0)$.
In the same way that one defines invariant closed subscheme of $(X,v)$, we define \textit{invariant foliations} by means of the Lie-derivative along the vector field $v$.

The invariant foliations of a smooth irreducible $D$-variety $(X,v)$ provide a effective tool to study rational factors --- namely rational dominant morphisms $\phi : (X,v) \dashrightarrow (Y,w)$ of $D$-varieties --- of the $D$-variety $(X,v)$. Indeed, we prove that any such rational factor induces through its tangent sheaf, an invariant foliation on $(X,v)$ (see Proposition \ref{invariant-tangentfoliation}).

We also prove that the singular locus of an invariant foliation is in fact invariant (as a closed invariant subvariety). We get this result by using the local complex analytic flow of the vector field $v$ instead of the vector field itself (see Corollary \ref{invariance-singular}).

When working dynamically in the fourth section, this convenient result will allow us to work far away from the singularities of the foliation. A similar invariance result, regarding the indeterminacy locus of a rational integral, was also needed in the proof of the criterion of orthogonality to the constants in \cite{moi}. 

\subsection{Definition}

{\defn\label{invariant-foliation} Let $X$ be an algebraic variety over $k$ and $v$ a vector field on $X$. We say that a coherent subsheaf $\mathcal F$ of $\Theta_{X/k}$ is \textit{$v$-invariant} if $\mathcal F$ is a $D$-coherent subsheaf of $(\Theta_{X/k},[v,-])$.}

In other words, a coherent subsheaf $\mathcal F$ of $\Theta_{X/k}$ is \textit{$v$-invariant} if
$$[v,\mathcal F(U)] \subset \mathcal F (U) \text{ for every open subspace }U \subset X.$$

{\rem\label{generic-invariance} Let $(\mathcal E,\nabla_\mathcal E)$ be a $D$-coherent sheaf over a $D$-scheme $(X,\delta_X)$. For every coherent subsheaf $\mathcal F \subset \mathcal E$, the derivation $\nabla_E$ induces a morphism of $\mathcal O_X$-sheaves:
$$\phi_\mathcal F:  \mathcal F \longrightarrow \mathcal E /\mathcal F.$$
Moreover, $\mathcal F$ is a $D$-coherent subsheaf of $\mathcal E$ if and only if this map is the null morphism.}

Now, consider the case of a linear differential equation defined over some constant differential field.

{\exam If $x_1, \ldots, x_n$ are coordinates on $\mathbb{A}^n$, then the vector fields $\frac \partial {\partial x_1}, \ldots, \frac \partial {\partial x_n}$   define a trivialization of the tangent bundle $T_{\mathbb{A}^n/k}$. We will work in this trivialization.

Consider a linear vector field $v(x) = A.x$ for some matrix $A \in \mathcal M_n(k)$
where $(k,0)$ is a constant differential field. Let $E \subset \mathbb{A}^n$ be a $A$-invariant linear subspace.
Then, the sheaf of sections $\mathcal E$ of $E \times \mathbb{A}^n$ is $v$-invariant as a subsheaf of $\Theta_{\mathbb{A}^n/k}$.

In particular, if $k$ is algebraically closed and the eigenvalues of $A$ are simple then the tangent bundle  
$$ \Theta_{\mathbb{A}^n/k} = \mathcal E_{\lambda_1} \oplus \cdots \oplus \mathcal E_{\lambda_n}.$$
splits as a direct sum of $1$-dimensional $v$-invariant coherent subsheaves.}

{\lem\label{generic-invariant-foliation} Let $X$ be an irreducible algebraic variety over $k$, let $\mathcal F \subset \Theta_{X/k}$ be a coherent saturated subsheaf and  let $U$ be a non-empty open subset of $X$.

The subsheaf $\mathcal F$ if $v$-invariant if and only if $\mathcal F_{|U}$ is $v_{|U}$-invariant.} 

\begin{proof}
The direct implication is obvious. Conversely, assume that $\mathcal F_{|U}$ is $v_{|U}$-invariant.
Using Remark \ref{generic-invariance}, we can consider the morphism of $\mathcal O_X$-coherent sheaves defined by the Lie-bracket with $v$:

$$ [v,-] : \mathcal F \longrightarrow \Theta_{X/k}/\mathcal F $$

By assumption, this morphism is the null morphism on $U$. Since $\Theta_{X/k}/\mathcal F$ does not have torsion, this implies that it is the the null morphism. 
\end{proof}
The principal examples of invariant foliations that we will encounter in those notes will be derived from rational factors by the following proposition.

{\Prop\label{invariant-tangentfoliation}  Let $\phi :(X,v)\dashrightarrow (Y,w)$ be a rational morphism of $D$-varieties over $(k,0)$. The tangent foliation $\mathcal F_{\phi}$ of $\phi$ is $v$-invariant.}

\begin{proof}
By Lemma \ref{generic-invariant-foliation}, we may assume that $X$ and $Y$ are smooth and that $\phi$ is regular and smooth. In particular, we have an exact sequence of locally free sheaves over $X$:
$$0 \longrightarrow \mathcal F_\phi \longrightarrow \Theta_{X/k} \overset{d\phi}{\longrightarrow} \phi^\ast \Theta_{Y/k} \longrightarrow 0$$

Moreover, since $\phi$ is a morphism of $D$-varieties, Lemma \ref{functoriality} ensures that $d\phi$ is a morphism of $D$-coherent sheaves over $(X,v)$:
$$ d\phi: (\Theta_{X/k}, [v,-]) \longrightarrow \phi^\ast (\Theta_{Y/k}, [w,-])$$

From Lemma \ref{abelian-category}, we conclude that its kernel $\mathcal F_\phi$ is a $D$-subcoherent sheaf of $\Theta_{X/k}$.
\end{proof}

\subsection{Analytification of an invariant coherent sheaf}

{\Prop\label{analytification invariance} Let $\mathcal F$ be a coherent analytic subsheaf of $\Theta_M$ and $v$ an analytic vector field on $M$. The following are equivalent:

\begin{itemize}
\item[(i)] The subsheaf $\mathcal F$ is $v$-invariant.
\item[(ii)] For every open set $U \subset M$ and $t$ such that the flow of $v$ is defined at $t$ for every $x \in U$,
$$\phi_t^\ast(\mathcal F(U)) \subset \mathcal F(\phi_{-t}(U)).$$
\end{itemize}}
 
 \begin{proof}
Using the semi-group properties of flows, the property (ii) is equivalent to:

(ii)' For every $x \in M$ and $t \in \C$ sufficiently small 
$$\phi_t^\ast(\mathcal F_x) \subset \mathcal F_{\phi_{-t}(x)}.$$

Denote by $\pi : \Theta_M \longrightarrow \Theta_M /\mathcal F$, the canonical projection.

Let $w \in \mathcal F_x$ be a local section of $\mathcal F$. Consider $\epsilon > 0$ such that $\phi^\ast_t w$ is defined and bounded on $B(x,\epsilon)$ for $t \in \mathbb C$ sufficiently small.

Now, the space of bounded sections of $\Theta_M$ and $\Theta_M /\mathcal F$ on $B(x,\epsilon)$ is a Banach space.

Using Cauchy formulas (Lemma \ref{Cauchyformulas}) we know that $t \mapsto \phi_t^\ast w$ is the solution of the differential equation: 
$$\begin{cases}
\phi_t^\ast(w)_{|t = 0} = w \\
\frac d {dt} \phi_t^\ast(w) = \mathcal L_v(\phi_t^\ast(w))
\end{cases}$$ 

Therefore, composing with the linear map $\pi$, the path $ t \mapsto \pi \circ \phi_t^\ast(w)$ is a solution of the differential equation:

$$\begin{cases}
\pi \circ \phi_t^\ast(w)_{|t = 0} = 0 \\
\frac d {dt} (\pi \circ \phi_t^\ast(w)) = \pi \circ \mathcal L_v(\phi_t^\ast(w))
\end{cases}$$ 

In the Banach space of bounded sections of $\Theta_M /\mathcal F$ on $B(x,\epsilon)$, we get that 

$$\pi \circ \phi^\ast_t w = 0 \text{ for } t  \text{ sufficiently small if and only if } \pi \circ \mathcal L_v(w) = 0$$
which exactly means the equivalence between (i) and (ii).
\end{proof}
 
{\cor\label{invariance-singular} Let $v$ be a vector field on $M$ and $\mathcal F$ an analytic foliation such that $\mathcal F$ is $v$-invariant.
Then, $\mathrm{Sing}(\mathcal F)$ is a closed invariant analytic subspace of $(M,v)$.}

\begin{proof}
By definition the singular locus of $\mathcal F$ is described by :
$$\mathrm{Sing}(\mathcal F) = \lbrace x \in M \text{ | }  \Theta_{M,x} / \mathcal F_x   \text{ is not a free }\mathcal O_{M,x}-\text{submodule of rank }1 \rbrace$$

For every point $p \in M$ and $t \in \mathbb{C}$ sufficiently small, the local flow $\phi_t$ at $p$ is a local diffeomorphism. It defines a ring-isomorphism between $\mathcal O_{X,x}$ and $\mathcal O_{X,\phi_t(x)}$ and its differential defines an isomorphism of modules: 
$$d\phi_t : \Theta_{M,x} \longrightarrow \Theta_{M,\phi_t(x)}$$

Since $\mathcal F$ is $v$-invariant,by Proposition \ref{analytification invariance}, we have that $\mathcal F_{\phi_t(x)} = \phi_t(\mathcal F_x)$. Therefore, $d\phi_t$ defines an isomorphism of exact sequences:
$$
    \xymatrix{0 \ar[r]  & \mathcal F_x \ar[r]\ar[d]  & \Theta_{M,x} \ar[r]\ar[d]^{d\phi_t} & \Theta_{M,x}/ \mathcal F_x \ar[r]\ar[d] & 0  \\
              0 \ar[r] & \mathcal F_{\phi_t(x)} \ar[r]  & \Theta_{M,\phi_t(x)} \ar[r] & \Theta_{M,\phi_t(x)}/ \mathcal F_{\phi_t(x)} \ar[r] & 0 }$$

Therefore, the first exact sequence is a direct product if and only if the second one is. We conclude that for $t$ sufficiently small, the germ of  $\mathrm{Sing}(\mathcal F)$ at $p$ is invariant by the local flow at $p$. Since this is true through any point of $\mathrm{Sing}(\mathcal F)$, we conclude that $\mathrm{Sing}(\mathcal F)$ is a closed invariant analytic subspace of $M$.
\end{proof}

{\cor[Preparatory Lemma]\label{preparation-corollary} Let $(X,v)$ be an absolutely irreducible real $D$-variety such that $M = X(\R)^{an}$ is regular, compact, and Zariski-dense in $X$.
Suppose that $(X,v)$ admits an invariant saturated coherent subsheaf $\mathcal F$ of $\Theta_{X/k}$. Then

\begin{itemize}
\item The singular locus  $Z = \mathrm{Sing}(\mathcal{F})$ is invariant.
\item On the dense invariant open subset $U = X(\R)^{an} \setminus Z(\R)$,  the real analytication of $\mathcal F$ defines a continuous vector subbundle $F$ of $T_{M}$  such that 
$$ \forall t \in \mathbb{R}\text{, } \forall x \in U \text{, } d\phi_t(F_x) = F_{\phi_t(x)}.$$
\end{itemize}}

\begin{proof}

Using Corollary \ref{invariance-singular}, the singular locus is invariant. Therefore, the restriction of the flow of $v$ to $U = X(\R)^{an} \setminus Z(\R)$  is complete. We then apply Proposition \ref{analytification invariance}.
\end{proof}

\section{Rational factors of mixing Anosov flows of dimension $3$}

If $(X,v)$ be a $D$-variety defined over $(\mathbb{R},0)$, we already reduced in the previous sections the understanding of the rational factors of $(X,v)$ to the understanding of the invariant algebraic foliations on $(X,v)$.

When the real-analytic flow $(M,(\phi_t)_{t \in \R})$ of $(X,v)$ is an Anosov flow, there exists, by definition, a splitting of the tangent bundle of $M$ into continuous invariant (non-singular) foliations: 
$$ T_M = W^{su} \oplus \mathbb{R}.v \oplus  W^{ss}$$  
where $W^{ss}$ (resp. $W^{su}$) is called the \textit{strongly stable} foliation (resp. \textit{strongly unstable} foliation).

In dimension $3$, these three invariant foliations have rank $1$ and therefore, can not be split again into invariant foliations of smaller rank. Looking at the periodic orbits of $(M,(\phi_t)_{t \in \R})$ , we are able to describe all continuous  invariant foliations of $(M,(\phi_t)_{t \in \R})$  from the strongly stable foliation $W^{ss}$, the strongly unstable one $W^{su}$ and the direction of the flow $\mathbb{R}.v$  (Proposition \ref{Anosov-invariantlinebundle} and Proposition \ref{Anosov-invariant2bundle}).

Using this explicit description of invariant continuous foliations and a result of Plante in \cite{Plante}, we conclude that none of these foliations comes from a rational factor.

\subsection{Continuous invariant subbundles of a mixing Anosov flow}

{\Prop \label{Anosov-invariantlinebundle} Let $(M,(\phi_t)_{t \in \R})$ be a compact and connected Anosov flow of dimension $3$ and $\Sigma$ a proper closed invariant subset. The $(d\phi_t)_{t \in \R}$-invariant continuous line subbundles defined on $U = M \setminus \Sigma$ are exactly:

\begin{itemize}
\item the (non-singular) foliation $F$ tangent to the flow.
\item the strongly stable and strongly unstable folations $W^{ss}$ and $W^{su}$, associated to the Anosov structure.
\end{itemize}}

\begin{proof}
By definition of an Anosov flow, these three continuous foliations are $(d\phi_t)_{t \in \mathbb{R}}$-invariant. We denote by $\sigma_1,\ldots, \sigma_2, \sigma_3 : M \longrightarrow \mathrm{Gr_1(TM)}$. the three (continuous) sections of the $1$-Grassmanian bundle of $TM$ defining those $1$-foliations.

Conversely let $L \subset TM$ be a continuous line subbundle defined on $U \subset M$ and denote by $\sigma : M \longrightarrow \mathrm{Gr_1(TM)}$, the associated section.

Since $(M,(\phi_t)_{t \in \mathbb{R}})$ is an Anosov flow, the periodic points of $(\phi_t)_{t \in \mathbb{R}}$ are dense in $U$.

Let $p \in M$  be a periodic points of $(\phi_t)_{t \in \mathbb{R}}$ of period $T > 0$. Then, the fibre $L_p \subset T_p M$ is a stable line of the linear map: 
$$(d\phi_{T})_{|T_pM} : T_p M \longrightarrow T_p M$$ 

Since $(M,(\phi_t)_{t \in \mathbb{R}})$ is an Anosov flow, this linear map has exactly three eigenvalues and the associated eigenspaces are precisely $W^{su}_p$, $W^{ss}_p$ and $F_p$.

We have proven that, on a dense set, $\sigma$ agrees with one of the sections $\sigma_1, \sigma_2, \sigma_3$. For $i \leq 3$, denote by $F_i$ be the closed  subset of $M$ where $\sigma$ agrees with $\sigma_i$.
Since $F_1 \cup F_2 \cup F_3$ is closed and contain a dense set, we have that $F_1 \cup F_2 \cup F_3 = M$. Moreover, they are disjoint since two distinct $\sigma_i$ have distinct values at every point of $M$. Since $M$ is connected, $M = F_i$ for some $i \leq 3$. This implies that $L$ is either the stable foliation, or the unstable foliation or the direction of the flow itself.
\end{proof}

{\Prop \label{Anosov-invariant2bundle} Let $(M,(\phi_t)_{t \in \R})$ be a compact and connected Anosov flow of dimension $3$ and $\Sigma$ a proper closed invariant subset. The $(d\phi_t)_{t \in \R}$-invariant continuous subbundles of rank $2$ defined on $U = M \setminus \Sigma$ are exactly:

\begin{itemize}
\item the stable and unstable folations $W^s = F \oplus W^{ss}$ and $W^u = F \oplus W^{su}$ associated to the Anosov structure, where $F$ is the direction tangent to the flow.
\item the direct sum of the strongly stable and strongly unstable folations $W^{ss} \oplus W^{su}$.
\end{itemize}}

\begin{proof}
Since these three examples are direct sum of  $(d\phi_t)_{t \in \mathbb{R}}$-invariant line bundles, there are also  $(d\phi_t)_{t \in \mathbb{R}}$-invariant.

Conversely let $P \subset TM$ be a continuous plane subbundle defined on $U \subset M$. Since $(M,(\phi_t)_{t \in \mathbb{R}})$ is an Anosov flow, the periodic points of $(\phi_t)_{t \in \mathbb{R}}$ are dense in $U$.

Let $p \in M$  be a periodic points of $(\phi_t)_{t \in \mathbb{R}}$ of period $T > 0$. Then, the fibre $L_p \subset T_p M$ is a stable plane of the linear map: 
$$(d\phi_{T})_{|T_pM} : T_p M \longrightarrow T_p M$$ 

But these stable planes are exactly $W^s_p$, $W^u_p$ and $(W^{ss} \oplus W^{su})_p$. One concludes in the same way as Proposition \ref{Anosov-invariantlinebundle}.
\end{proof}

{\Thm[{\citep[Theorem 1.3]{Plante}}]\label{Plante} Let $(M,(\phi_t)_{t \in \R})$ be a mixing Anosov flow.

Every leaf of the strongly stable foliation $W^{ss}$ and every leaf of the strongly unstable foliation $W^{su}$ is dense in $M$.}

We will use the Theorem \ref{Plante} combined with the two preceding propositions in the following form:

{\cor\label{key lemma} Let $(M,(\phi_t)_{t \in \R})$ be a compact and connected mixing Anosov flow of dimension $3$, $\Sigma$ a proper closed invariant subset and $F$ a continuous foliation on $M \setminus \Sigma$ with positive rank.

Suppose that the foliation $F$ is $(d\phi_t)_{t \in \R}$-invariant and distinct from the foliation tangent to the flow. Then, every leaf of $F$ is dense in $M$.}

\begin{proof}
Denote by $r$ the rank of the foliation. If $r = 3$, it is true since $\Sigma$ has empty interior (the flow is mixing so topologically transitive). For $r = 1,2$, we do a case by case inspection:
\begin{itemize}
\item If $r = 1$, then by Proposition \ref{Anosov-invariantlinebundle}, $F$ is either the  strongly stable foliation or the strongly unstable one. In both cases, Theorem \ref{Plante} implies that every leaf of $F$ is dense.

\item If $r = 2$, then by Proposition \ref{Anosov-invariant2bundle}, $F$ is either the stable foliation or the strongly unstable one, or the direct sum of the strongly stable and the strongly unstable ones. In those three cases, the foliation contains either the strongly stable foliation or the strongly unstable one. Using Theorem \ref{Plante}, we conclude that every leaf of $\mathcal F$ is dense in $M$.

\end{itemize}
 
\end{proof}
\subsection{Rational factors of mixing, compact, Anosov flows of dimension $3$} 

{\Thm\label{maintheoremcore} Let $(X,v)$ be an absolutely irreducible $D$-variety of dimension $3$ over $\mathbb{R}$. Assume that the real-analytification $X(\R)^{an}$ of $X$ admits a compact (non-empty) connected component $C_\mathbb{R}$ contained in the regular locus of $X$.

If the real analytic flow $(C_\mathbb{R}, (\phi_{t |C_\mathbb{R}})_{t \in \R})$ is a mixing Anosov flow, then $(X,v)$ does not admit any non-trivial rational factor.}

\begin{proof}
Let $(X,v)$ be a real $D$-variety satisfying the hypothesis of Theorem \ref{maintheoremcore}. Suppose that $(X,v)$ admits a non-trivial rational factor $\pi : (X,v) \dashrightarrow (Y,w)$. Since we already now that $(X,v)$ has no non-trivial rational integral (see \cite{moi2}), we may assume that $w \neq 0$. By Proposition \ref{invariant-tangentfoliation}, the tangent foliation $\mathcal F_{\pi}$ is $v$-invariant and does not contain the foliation generated by $v$, since $w \neq 0$. 

Now, by corollary \ref{singularlocus}, the singular locus $Z$ of $\mathcal F_{\pi}$ is a closed invariant subvariety. On the open set $U = X \setminus Z$, $\mathcal F_{\pi}$ is a non-singular foliation. Consequently, outside of the closed invariant set $\Sigma = Z(\R)$, $F = \mathcal F_{\pi}^{an}$ is a non-singular foliation on $X(\R)^{an} \setminus \Sigma$.

\begin{itemize}
\item Since $\mathcal F_{\pi}$ is $v$-invariant, the continuous foliation is a $(d\phi_t)_{t \in \R}$-invariant subbundle of $T_{M}$.
\item It is distinct from the foliation generated by $v$, since $w \neq 0$.
\end{itemize}

By Corollary \ref{key lemma}, we conclude that every leaf of $F$ is dense in $C_\R$. This contradicts the fact that $\mathcal F_{\pi}$ is a foliation tangent to the algebraic fibration $\pi$ (see Example \ref{example}).
\end{proof}

\subsection{Proof of Theorem \ref{maintheoremInt}} 

{\Thm\label{final} Let $X$ be an absolutely irreducible variety of dimension $3$ over $\mathbb{R}$ endowed a vector field $v$. Assume that the real-analytification $X(\R)^{an}$ of $X$ admits a compact (non-empty) connected component $C_\mathbb{R}$ contained in the regular locus of $X$.

If the real-analytic flow $(C_\mathbb{R}, (\phi_t)_{t \in \R})$ is a mixing Anosov flow, then the autonomous differential equation $(X,v)$ is generically disintegrated.}

\begin{proof} 
Consider an algebraic autonomous differential equation $(X,v)$ satisfying the hypotheses of Theorem \ref{final}. It is easy to see that the algebraic differential equation $(X,v)$ also satisfies the hypotheses of  Théorème D of \cite{moi}: the differential equation $(X,v)$ is absolutely irreducible and the dynamic of the real-analytic flow of the vector field $v$ on the invariant compact subspace $K = C_\R$ is mixing (hence, weakly mixing).

It therefore follows from Théorème D that the generic type of $(X,v)$ is orthogonal to the constants and therefore from Théorème A of \cite{moi} that the differential equation $(X,v)$ admits a rational factor $\pi: (X,v) \dashrightarrow (Y,w)$ of positive dimension, which is generically disintegrated.

Now, since $Y$ has positive dimension, Theorem \ref{maintheoremcore} ensures that the generic fibre of $\pi$ is in fact finite, so that $(X,v)$ is a finite extension of the generically disintegrated differential equation $(Y,w)$. We conclude using Proposition 1.3.7 of \cite{moi} that $(X,v)$ itself is generically disintegrated.
\end{proof}

\subsection{Proof of Corollary \ref{maintheoremInt}} 

{\cor\label{maincorollarycore} Let $M$ be a regular compact real-algebraic subset of the Euclidean space $\mathbb{R}^N$ of dimension $2$ with negative curvature and let $r$ be integer $ \geq 2$.

Consider $r$ unitary geodesics $\gamma_1, \ldots \gamma_r: \mathbb{R} \rightarrow SM$ of  the Euclidean submanifold $M$, viewed as analytic curves on the sphere bundle $SM \subset T \mathbb{R}^N$ of $M$. Assume that, for every $1 \leq i \leq r$, the analytic curve $\gamma_i$ is Zariski-dense in $SM$. Then, the following are equivalent:
\begin{itemize}
\item[(i)] The analytic curve $t \mapsto (\gamma_1(t), \ldots, \gamma_r(t))$ is Zariski-dense in $SM^r$.
\item[(ii)] For every $i \neq j$, the analytic curve $t \mapsto (\gamma_i(t),\gamma_j(t))$ is Zariski-dense in $SM^2$.  
\end{itemize}}

Note that $(i) \Longrightarrow (ii)$ is always true and that the real content of Corollary \ref{maincorollarycore} is contained in the converse implication.

\begin{proof}

Consider a regular compact real-algebraic subset  $M$ of the Euclidean space $\mathbb{R}^N$ satisfying the hypotheses of Corollary \ref{maincorollarycore}. We write $M = X(\R)$ for some quasi-affine algebraic $X$ over $\mathbb{R}$.

Let's first argue that we may assume that $X$ is absolutely irreducible: since $M$ is regular and connected, it is contained in an irreducible component of $X$, so that we may assume that $X$ is irreducible. Moreover, by replacing $X$ by the Zariski-closure of $M$ in $X$, we can assume that $M = X(\R)$ is Zariski-dense in $X$. Together with the irreducibility property, this implies that $X$ is absolutely irreducible (see for example Lemma 4.1.1 of \cite{moi2}). 

The Euclidean metric on the affine space endows $X$ with the structure of a pseudo-Riemmannian manifold over $\R$ in the sense of \cite{moi2}. Denote by $(SX,v)$ the corresponding (unitary) geodesic differential equation (supported on the sphere bundle $SX$ of $X$).

{\asser With the notation above, the autonomous differential equation $(SX,v)$ is generically disintegrated.}

\begin{proof} Indeed, it follows from the proof of Theorem 4.2.1 in \cite{moi2} that:
\begin{itemize}
\item The subset $SX(\mathbb{R})$ is Zariski dense in $SX$.
\item The algebraic variety $SX$ (over $\R$) is absolutely irreducible.
\item The real-analytic flow of $(SX,v)$ is the unitary geodesic flow of the compact Riemanian manifold $(M,g_{E|M})$ where $g_{E|M}$ denotes the restriction of the Euclidean metric to $M$.
\end{itemize}

Since $(M,g_{E|M})$ is a compact Riemannian manifold of dimension two with negative curvature, it follows from the results of \cite{Ano} that the unitary geodesic flow of $(SX,v)$ is a mixing Anosov flow of dimension $3$ (see, for example, Theorem 3.3.5 of \cite{moi2}). We conclude that $(X,v)$ satisfies the hypotheses of Theorem \ref{maincorollarycore} and is therefore generically disintegrated.
\end{proof}

Now consider $\gamma_1, \ldots \gamma_r: \mathbb{R} \longrightarrow SM$, geodesics of $M$, Zariski-dense in $SM$, such that the analytic curve $\Gamma: t \mapsto (\gamma_1(t), \ldots, \gamma_r(t))$ is not Zariski-dense in $SM^r$ and denote by $G$ its Zariski closure is $SX^r$. The following properties are fulfilled:

\begin{itemize}
\item $G$ is a proper closed algebraic subvariety of $SX^r$ which projects generically on all factors.

\item $G$ is an invariant closed subvariety of $(SX,v)^r$: indeed, since $\Gamma$ is invariant under the action of the real-analytic flow of the vector field $v \times \ldots \times v$, the Zariski-closure $G$ of $\Gamma$ is a proper closed invariant subvarieties of $(SX,v)^r$. 
\end{itemize}

Since $(SX,v)$ is generically disintegrated, we can write $G$ as an irreducible component of:
$$\bigcap_{1 \leq i \neq j \leq r} \pi_{i,j}^{-1}(Z_{i,j}).$$
where $\pi_{i,j} : X^r \longrightarrow X^2$ is the projection on the $i^{th}$ and $j^{th}$ coordinates and $Z_{i,j} \in \mathcal I^{gen}_2(SX,v)$ for every $i \neq j \leq r$.

Since $G$ is a proper closed subvariety of $(SX,v)^r$, for some $i \neq j \leq n$, the closed subvariety  $Z_{i,j}$ is a proper closed subvariety of $(SX,v)^2$. It follows that $t \mapsto (\gamma_i(t),\gamma_j(t))$ is not Zariski-dense in $SM^2$, which concludes the proof.
\end{proof}
\bibliographystyle{alpha}
\bibliography{bibliographie}
\end{document}